%% file: main.tex
\documentclass[conference]{IEEEtran}
\usepackage{bm}
\usepackage{epsfig}
\usepackage{graphicx}
\usepackage{subfigure}
\usepackage{float}
\usepackage{cite}
\usepackage{url}
\usepackage{color}
\usepackage{balance}
\usepackage{mdwlist}
\usepackage{multirow}
\usepackage{threeparttable}
\usepackage{enumitem}
\usepackage{amsmath}
\usepackage{amssymb,lipsum}
\usepackage{stmaryrd}
\usepackage{booktabs}
\usepackage{cuted}
\usepackage{siunitx}
\usepackage[numbers,sort&compress]{natbib}
\usepackage{threeparttable}
\usepackage{epsfig,amsmath,amsfonts}
\usepackage{amsfonts}

\newcommand{\hp}{\mathop{\mathrm{hp}}}

\def\R{\mathbb{R}}

\makeatletter
\newif\if@restonecol
\makeatother

\usepackage[ruled,vlined]{algorithm2e}
\usepackage{algpseudocode}
\usepackage{amsmath}

\newcommand{\parm}{{\xi}}
\newcommand{\vecpar}{\boldsymbol{\parm}}

\newcommand{\multiGPC}{\Psi }

\newcommand{\polyInd}{\alpha}
\newcommand{\basisInd}{\boldsymbol{\polyInd}}

\newcommand{\ten}[1]{\mathcal{#1}}

\newcommand{\mat}[1]{\mathbf{#1}}

\newcommand{\reff}[1]{(\ref{#1})}

\RequirePackage{url}
\def\BibTeX{{\rm B\kern-.05em{\sc i\kern-.025em b}\kern-.08em
    T\kern-.1667em\lower.7ex\hbox{E}\kern-.125emX}}
\begin{document}

\IEEEoverridecommandlockouts
\title{High-Dimensional Uncertainty Quantification via Active and Rank-Adaptive Tensor Regression
}
\IEEEpubidadjcol
\author{\IEEEauthorblockN{Zichang He and Zheng Zhang}
\IEEEauthorblockA{Department of Electrical and Computer Engineering, 
University of California, Santa Barbara, CA 93106\\
Emails: zichanghe@ucsb.edu, zhengzhang@ece.ucsb.edu}
}


\maketitle
\begin{abstract}
Uncertainty quantification based on stochastic spectral methods suffers from the curse of dimensionality. This issue was mitigated recently by low-rank tensor methods. 
However, there exist two fundamental challenges in low-rank tensor-based uncertainty quantification: how to automatically determine the tensor rank and how to pick the simulation samples. This paper proposes a novel tensor regression method to address these two challenges. Our method uses a $\ell_{q}/ \ell_{2}$-norm regularization to determine the tensor rank and an estimated Voronoi diagram to pick informative samples for simulation. 
The proposed framework is verified by a 19-dim photonic band-pass filter and a 57-dim CMOS ring oscillator, capturing the high-dimensional uncertainty well with only 90 and 290 samples respectively.
\end{abstract}


\input{introduction.tex}

\input{preliminary.tex}
\input{proposed_method.tex}
\input{numerical_result.tex}
\input{conclusion.tex}



\small{
\bibliographystyle{IEEEtran}
\bibliography{Bib/reference}
}

\end{document}

%% file: introduction.tex
\section{Introduction}\label{sec:introduction}
Fabrication process variations are a major concern in nano-scale chip design. To estimate and quantify the uncertainties caused by process variations, Monte Carlo (MC) is the mainstream uncertainty quantification (UQ) tool used in commercial EDA tools, but it requires a huge amount of simulation samples. Instead, stochastic spectral methods based on generalized polynomial chaos (gPC)~\cite{xiu2010numerical} offer an efficient alternative by approximating a stochastic circuit performance metric as a linear combination of some basis functions~\cite{zhang2013stochastic,He2019,ginste2012stochastic}. However, stochastic spectral methods suffer from the curse of dimensionality: a huge amount of simulation samples are required when the number of random parameters is large. 

Low-rank tensor methods are a promising technique to solve high-dimensional UQ problems~\cite{konakli2016polynomial,zhang2014enabling, zhang2016big,chevreuil2015least}. In~\cite{zhang2014enabling}, a high-dimensional gPC expansion is obtained via a low-rank tensor recovery, which estimates massive unknown output samples from a few simulation results. However, the method~\cite{zhang2014enabling} uses a fixed tensor rank, which is hard to estimate {\it a-priori} in practice. The most recent work~\cite{shi2019meta} uses a greedy rank-1 update until a good accuracy is reached. However, greedy rank-1 tensor update does not provide optimal solutions and can cause over fitting. Besides, it is not clear how to adaptively pick the simulation samples to reduce the computation budget.

{\bf Contributions.} This paper proposes a novel high-dimensional UQ solver based on tensor regression. 
In order to automatically determine the tensor rank, we employ a group-sparsity regularization in the training process.  
We also develop an adaptive sampling strategy to reduce the simulation cost. This method balances exploration and exploitation of our model. Our method is used to quantify the uncertainties of a 19-dim phonic IC and a 57-dim electronic IC with 90 and 290 simulation samples respectively.

%% file: preliminary.tex
\section{Background}
\label{sec:preliminaries}

{\bf Generalized Polynomial Chaos.} Let $\vecpar = \left[\parm_1,\ldots,\parm_d \right] \in \R^d$ be a random vector describing process variations. We aim to estimate the interested performance metric $y(\vecpar)$ (e.g., chip frequency or power) under such uncertainty.
A truncated gPC expansion approximates $y(\parm)$ as
\begin{equation}
\label{eq:Multi_expansion}
    y(\vecpar)\approx \hat{y}(\vecpar) = \sum_{\basisInd \in \Theta} {c}_{\basisInd}\multiGPC_{\basisInd}(\vecpar),
\end{equation}
where $\basisInd$ is an index vector, and $\multiGPC_{\basisInd}$ is a polynomial basis function of degree $|\basisInd|=\alpha_1 + \alpha_2+ \cdots +\alpha_d$. If the joint probability density function of $\vecpar$ is $\rho(\vecpar)$, then the basis functions satisfy the orthornormal condition:
\begin{equation}
    \langle \multiGPC_{\basisInd}(\vecpar), \multiGPC_{\boldsymbol{\beta}}(\vecpar) \rangle= \int \limits_{\mathbb{R}^d} \multiGPC_{\basisInd}(\vecpar) \multiGPC_{\boldsymbol{\beta}}(\vecpar) \rho(\vecpar)d \vecpar=\delta_{\basisInd,\boldsymbol{\beta}}.
\end{equation}
Once the index set $\Theta$ is chosen, we need to determine the unknown coefficient $c_{\basisInd}$ for each $\basisInd \in \Theta$. The gPC only requires a small number of basis functions and simulation samples when the parameter dimensionality $d$ is small. However, a huge number of basis functions and simulation samples are required when $d$ is large. For instance, in the classical stochastic collocation method~\cite{xiu2010numerical}, the  number of simulation samples required to obtain $c_{\basisInd}$'s is an exponential or polynomial function of $d$.

{\bf Tensors.} A promising tool to overcome the curse of dimensionality is the tensor. A $d$-dim tensor $\ten{X} \in \R^{n_1 \times \cdots n_d}$ represents a $d$-dimensional data array, and it becomes a matrix when $d=2$. The $(i_1,\cdots, i_d)$-th element of $\ten{X}$ can be denoted as $x_{i_1 \cdots i_d}$. Given two tensors $\ten{X}$ and $\ten{Y}$ of the same size, their inner product is defined as:
\begin{equation}\label{eq:inner_pro}
    \langle \ten{X}, \ten{Y} \rangle := \sum_{i_1\cdots i_d} x_{i_1\cdots i_d} y_{i_1\cdots i_d}.
\end{equation}
A $d$-dim rank-$R$ tensor can be written as the sum of $R$ rank-1 tensors, known as a CP decomposition:
\begin{equation}
\label{eq:CP}
    \ten{X}=\sum_{r=1}^R \mat{u}_r^{(1)} \circ \mat{u}_r^{(2)} \cdots \circ \mat{u}_r^{(d)} = [\![\mat{U}^{(1)}, \mat{U}^{(2)}, \ldots, \mat{U}^{(d)}]\!],
\end{equation}
where $\circ $ denotes an outer product. The last term is the Krusal form, where factor matrix $\mat{U}^{(k)} = \left[\mat{u}_1^{(k)}, \ldots, \mat{u}_R^{(k)} \right] \in \mathbb{R}^{n_k \times R}$ includes all vectors associated with the $k$-th dimension. 



%% file: proposed_method.tex
\section{Proposed Tensor Regression Method}
\label{sec:proposed_method}

\subsection{Tensor Regression Formulation}
We choose the following index set for the gPC expansion:
\begin{equation}
\label{eq:indesSet}
    \Theta=\left\{ \basisInd=[\alpha_1, \alpha_2, \cdots, \alpha_d]\; |\; 0\leq \{\alpha_i\}_{i=1}^{d} \leq p \right\}.
\end{equation}
This specifies a gPC expansion with $(p+1)^d$ basis functions. Let $i_k=\alpha_k+1$, then we can define two $d$-dimensional tensors $\ten{X}$ and $\ten{B}(\vecpar)$ with their $(i_1, i_2, \cdots i_d)$-th element specified as
\begin{align}
\label{eq:tensorDef}
x_{i_1 i_2 \cdots i_d}=c_{\basisInd} \; {\text{and}}\; b_{i_1 i_2 \cdots i_d}(\vecpar)=\multiGPC_{\basisInd}(\vecpar).
\end{align}

Combining \eqref{eq:Multi_expansion}, \eqref{eq:indesSet} and \eqref{eq:tensorDef}, the truncated gPC expansion can be written as a tensor regression model
\begin{equation}
   y(\vecpar) \approx \hat{y}(\vecpar) = \langle \ten{X}, \ten{B}(\vecpar) \rangle.
\end{equation}
It is worth noting that tensor $\ten{B}(\vecpar)$ depends on $\vecpar$. When the random parameters $\vecpar$ are mutually independent, $\multiGPC_{\basisInd}(\vecpar)$ can be written as the product of $d$ uni-variable basis functions for each parameter $\xi_k$. In this case $\ten{B}(\vecpar) $ is a rank-1 tensor.

Our goal is to compute $\ten{X}$ given a set of data samples $\left \{\vecpar_n, y(\vecpar_n)\right \}_{n=1}^N$. Assume that $\ten{X}$ has the rank-$R$ decomposition in (\ref{eq:CP}), we can solve the following optimization problem
\begin{equation}
\label{eq:coe_ten}
    \min \limits_{\{ \mat{U}^{(k)}\}_{k=1}^d} h(\ten{X})= \frac{1}{2} \sum \limits_{n=1}^N \left( y_n-\langle [\![\mat{U}^{(1)}, \mat{U}^{(2)}, \ldots, \mat{U}^{(d)}]\!], \ten{B}_n \rangle \right )^2,
\end{equation}
where $y_n=y(\vecpar_n)$ and $\ten{B}_n=\ten{B}(\vecpar_n)$.



\subsection{Automatic Rank Determination}
The low-rank tensor regression~\reff{eq:coe_ten} requires the rank of $\ten{X}$ to be determined {\it a-priori}, which is often infeasible in practice. In order to address this issue, we first choose a sufficiently large $R$ such that it is above the actual rank, then we choose a proper rank-shrinking penalty function to regularize~\reff{eq:coe_ten}.

Specifically, we employ a group $\ell_{q}/ \ell_{2}$-norm regularization function to shrink the rank of $\ten{X}$: 
\begin{equation}
\label{eq:rank_penalty}
   g(\ten{X}) = \|\mat{v}\|_q, 
   \mat{v}= \left( \sum \limits_{k=1}^d \| \mat{u}_{r}^{(k)} \|_2^2 \right)^{\frac{1}{2}} , \quad q \in \left( 0, 1\right].
\end{equation}
This function puts $\mat{u}_r^{(k)}$, the $r$-th column of each $\mat{U}^{(k)}$, in the same group, and measures the $\ell_{q}/ \ell_{2}$ norm of all groups. As a result, one can shrink some groups to zero by reducing $g(\ten{X})$, leading to an automatic rank reduction. A smaller $q$ leads to a stronger shrinkage, and $q=1$ corresponds to a group lasso. 

By adding the penalty term \reff{eq:rank_penalty}, we have the following improved tensor regression model:
\begin{equation}
\label{eq:tensor_regression}
\begin{aligned}
        \min \limits_{\{\mat{U}^{(k)}\}_{k=1}^d} f(\ten{X}) = & h(\ten{X}) + \lambda g(\ten{X}).
\end{aligned}
\end{equation}
After solving this optimization problem, each obtained factor matrix $\mat{U}^{(k)}$ has a few common columns whose values are close to zero. These columns can be deleted and the actual rank of our obtained tensor becomes $\hat{R}\leq R$, where $\hat{R}$ is the number of remaining columns that are not deleted.

It is non-trivial to minimize $f(\ten{X})$
since the regularization function $g(\ten{X})$ is usually non-differentiable and non-convex with respect to $\mat{U}^{(k)}$'s. 
Instead, we solve the following optimization problem in practice:
\begin{equation}
    \label{eq:tensor_regression_bound}
\begin{aligned}
        \min \limits_{\{\mat{U}^{(k)}\}_{k=1}^d} \hat{f}(\ten{X}) = & h(\ten{X}) + \lambda \hat{g}(\ten{X}).
\end{aligned}
\end{equation}
Here $\hat{g}(\ten{X})$ is a transformation of $g(\ten{X})$ obtained via the variational equality~\cite{jenatton2010structured}:
\begin{align}
\label{eq:subproblem_eta}
     \hat{g}(\ten{X}) = g(\ten{X}) = \min_{\boldsymbol{\eta} \in \R_{+}^R} \frac{1}{2}\sum \limits_{r=1}^R \frac{\sum \limits_{k=1}^d \| \mat{u}_{r}^{(k)} \|_2^2 }{{\eta}_r} + \frac{1}{2}\|\boldsymbol{\eta}\|_{\frac{q}{2-q}}.
\end{align}
Once $\{\mat{U}^{(k)}\}_{k=1}^d$ is given, the values of $\hat{g}(\ten{X})$ and $\hat{f}(\ten{X})$ can be estimated by setting the elements of $\boldsymbol{\eta}$ as
\begin{equation}
\label{eq:eta_solution}
    {\eta}_r = ({{z}_r})^{2-q}\| \mat{z} \|_{q}^{q-1}, \quad \forall \; r=1,\ldots,R,
\end{equation}
where ${z}_r = {\left( \sum \limits_{k=1}^d \| \mat{u}_{r}^{(k)} \|_2^2\right)}^{\frac{1}{2}}$.   

Problem~\reff{eq:tensor_regression_bound} can be solved via an alternating algorithm such as a block coordinate descent solver or alternating direction method of multipliers. Due to the page limitation, we omit the details in this paper and will explain the detailed optimization algorithm in an extended journal paper.

\subsection{Adaptive Sampling Strategy}
Another fundamental question is how to select the parameter samples $\{ \vecpar_n\}_{n=1}^N$ for simulation. The method in~\cite{zhang2016big} uses some Monte Carlo random samples. Instead, this paper reduces the simulation cost by selecting only a few informative samples for detailed device- or circuit-level simulations.

We first use the Latin Hybercube (LH) sampling to generate 
an initial sample set ${\Omega}$. Then we employ an exploration step via the Voronoi diagram to measure the sample density in $\Omega$.
Given two distinct samples $\vecpar_i, \vecpar_j \in \Omega$, a Voronoi cell $C_i(\vecpar_i)$ covers the region that are closest to $\vecpar_i$. It is defined as the intersection of a set of half-planes (hp):
\begin{equation}
\begin{split} 
         C_i(\vecpar_i) = &\bigcap_{\vecpar_j\in \Omega \setminus \vecpar_i} \hp(\vecpar_i,\vecpar_j)\\
         \hp(\vecpar_i,\vecpar_j)=&\{\vecpar \in \R^d | \; \|\vecpar-\vecpar_i\| \le \|\vecpar-\vecpar_j\|\}.
\end{split}
\end{equation}    
It is intractable to calculate an Voronoi cell exactly in a high-dimensional space. 
However, we can easily estimate it via Monte Carlo~\cite{crombecq2009space}. 
The sample density of $C_i$ is approximated by counting the number of samples that are closest to $\vecpar_i$. 
Each sample in $\Omega $ determines one Voronoi cell with itself as the center, and we can select a new sample from the cell region with the lowest density.

If the performance metric $y(\vecpar)$ is known to be highly nonlinear, we can further exploit its non-linearity. Given $\vecpar$ and a Voronoi cell center $\mat{a}$, we measure the nonlineary of $y(\vecpar)$ as
\begin{equation}
    \gamma(\vecpar) = |\hat{y}(\vecpar) - \hat{y}(\mat{a}) - \nabla \hat{y}(\mat{a})^T(\vecpar-\mat{a}) |.
\end{equation}
We select a new sample as the one with largest $\gamma(\vecpar)$ in a Voronoi cell with the lowest sample density.
This method can be easily extended to a batch version by searching the top-$K$ least-sampled regions.

%% file: numerical_result.tex
\section{Numerical Experiment}
\label{sec:numerical_results}
\begin{figure}[t]
    \centering
    \includegraphics[scale=0.3]{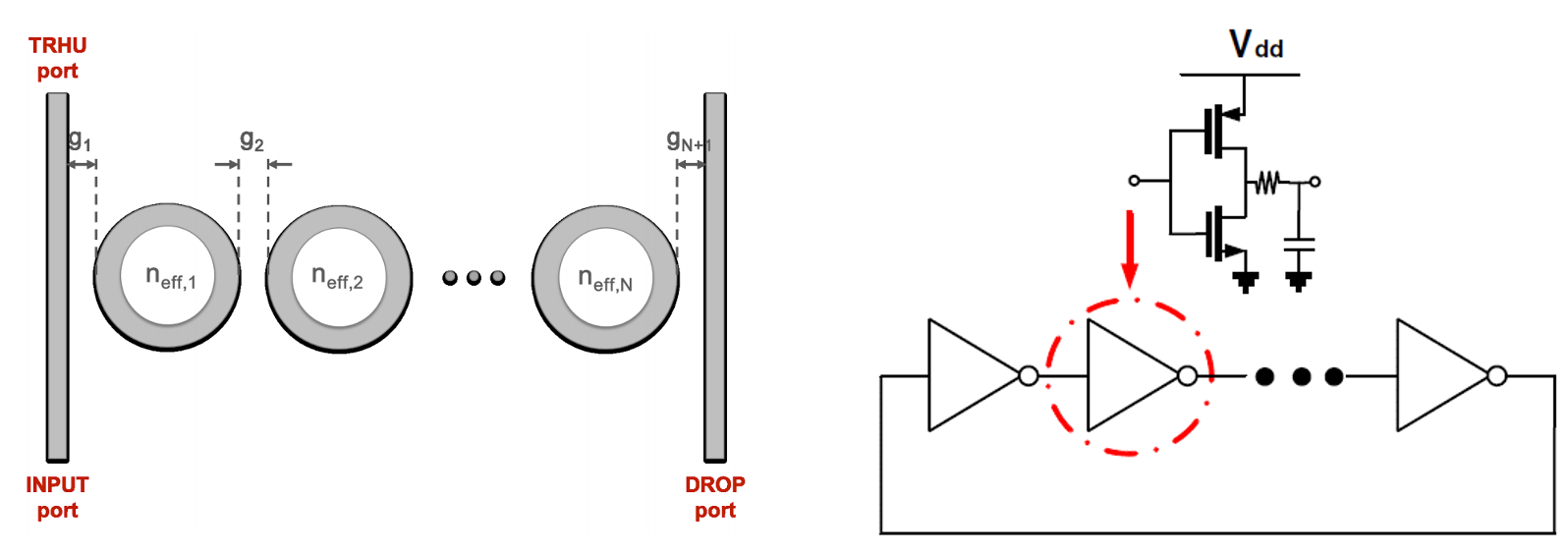}
    \caption{Left: a photonic band-pass filter with 9 micro-ring resonators. Right: Schematic of a CMOS ring oscillator.}
    \label{fig:PIC_CMOS}
\end{figure}

We verify our algorithm by a photonic band-pass filter and a CMOS ring oscillator as shown in Fig.~\ref{fig:PIC_CMOS}. 

\begin{figure}[t]
    \centering
    \includegraphics[width=1\columnwidth]{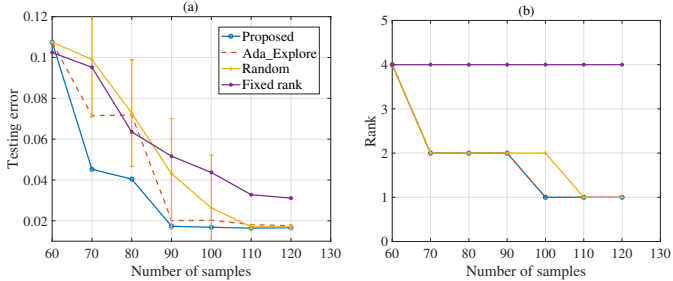}
    \caption{Result of the photonic filter. (a) Testing error on $10^5$ MC samples, achieving 98.3\% accuracy. (b) The estimated tensor ranks.}
    \label{fig:pic19}
\end{figure}

\textbf{Photonic band-pass filter:} This benchmark has 19 Gaussian random parameters describing the variations of the effective phase index ($n_\text{neff}$) of each ring, and the gap ($g$) between adjacent rings and between the end ring and the bus waveguides. 
We set the highest polynomial order $p=2$ for each dimension, and initialize $\ten{X}$ as rank-4, and use $q=0.5$ in the regularizer. 
We initialize 60 samples with a standard LH experimental design, and further select 6 batches with 60 additional samples. The proposed tensor regression framework is compared with a fixed rank method, a random and an adaptive exploration sampling method (shown in Fig.~\ref{fig:pic19}).

\begin{figure}[t]
    \centering
    \includegraphics[width=1\columnwidth]{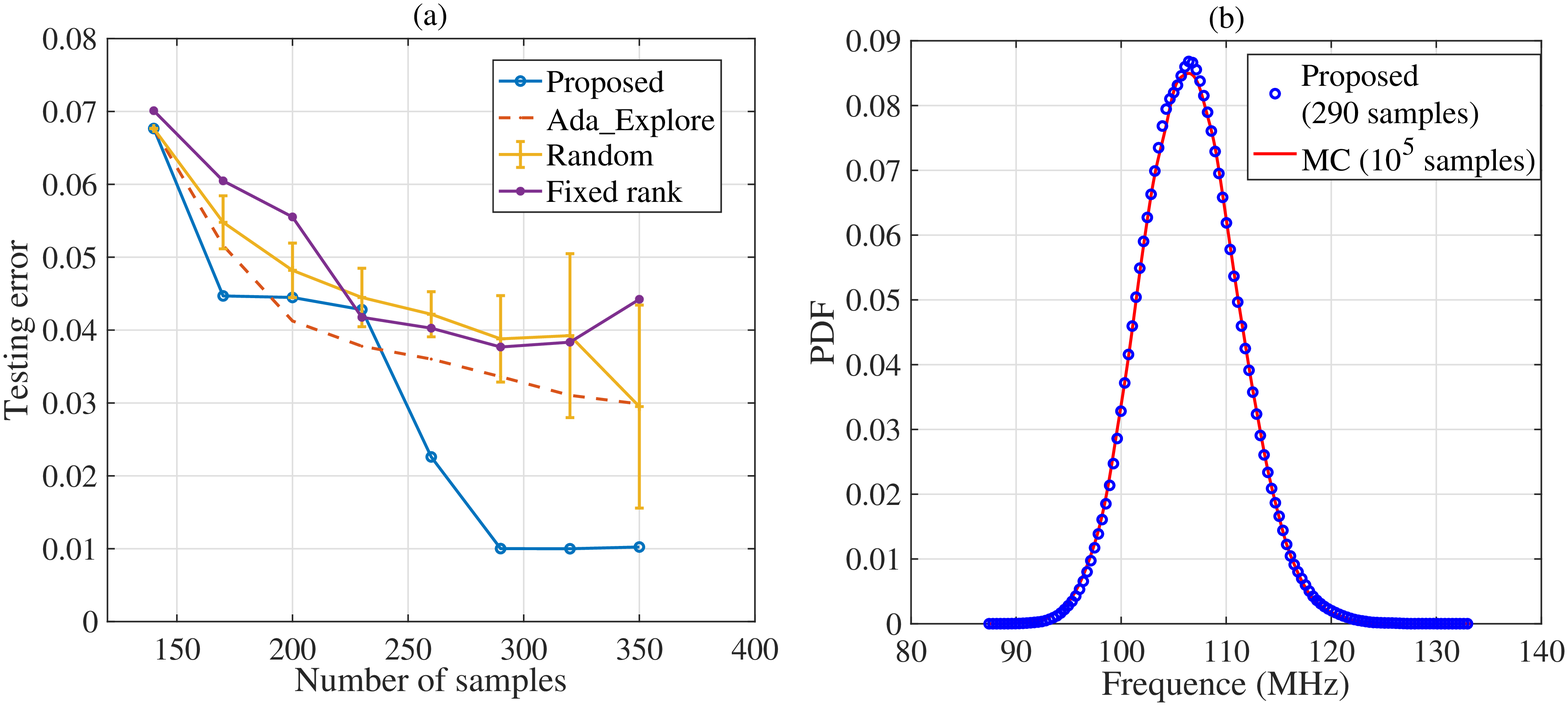}
    \caption{Results of the CMOS ring oscillator. (a) Testing error on $10^5$ MC samples (b) Probability density functions of the oscillator frequency.}
    \label{fig:adp57}
\end{figure}

\textbf{CMOS ring oscillator:} This circuit has 57 Gaussian random parameters describing the variations of threshold voltages, gate-oxide thickness, and effective gate length/width. 
The simulation results are obtained by calling a periodic steady-state simulator repeatedly.
We set basis as order-2 in each dimension, initialize $\ten{X}$ as rank-5, and use $q=0.5$ in the regularizer.
We start from 140 standard LH samples, and adaptively select additional 7 batches with 210 samples. 
The results and comparison are summarized in Fig.~\ref{fig:adp57} and Table~\ref{tab:ring57}. 

\begin{table}[t]
\centering
\caption{Comparison with the gPC model with a total degree scheme.}
\label{tab:ring57}
\begin{tabular}{ccccc}
\toprule
& Proposed  & Total-degree gPC  & MC \\
\midrule
\# of variables & 855     & 1711         &  N/A      \\
\# of samples   & 290       & 1711         &  $10^5$   \\
    Mean        & 106.28  & 106.58    &  106.53 \\
  Stdvar        & 4.616     & 6.810      &  4.641     \\
   Error        & 1\%       & 4.84\%        &  N/A      \\
\bottomrule
\end{tabular}
\vspace{-10pt}
\end{table}

%% file: conclusion.tex
\section{Conclusion}
This paper has proposed a tensor regression framework for high-dimensional uncertainty quantification. Our method has addressed two fundamental challenges: automatic tensor rank determination and adaptive sampling. 
The numerical result has demonstrated the excellent capability of automatic rank determination of our method, and the simulation cost reduction by our adaptive sampling method.
